\documentclass[11pt]{article}
\usepackage{amssymb,amsfonts,amsmath,latexsym,epsf,tikz,url}

\newtheorem{theorem}{Theorem}[section]

\newtheorem{lemma}[theorem]{Lemma}

\newcommand{\proof}{\noindent{\bf Proof.\ }}
\newcommand{\qed}{\hfill $\square$\medskip}

\begin{document}

\title{Domination polynomial of lexicographic product of specific graphs}

\author{
Saeid Alikhani $^{}$\footnote{Corresponding author}
\and
Somayeh Jahari
}

\date{}

\maketitle

\begin{center}
 Department of Mathematics, Yazd University, 89195-741, Yazd, Iran\\
{\tt alikhani@yazd.ac.ir}\\

\end{center}


\begin{abstract}
Let $G$ be a simple graph of order $n$.  The domination polynomial of $G$ is the polynomial $D(G,\lambda)=\sum_{i=0}^{n}
d(G,i) \lambda^{i}$, where $d(G,i)$ is the number of dominating
sets  of $G$ of size $i$.  We consider the lexicographic product of two specific graphs and study
  their domination polynomials.
\end{abstract}


\medskip
\noindent{\bf AMS Subj. Class.:} 05C60.

\section{Introduction}

 Let $G=(V, E)$ be a simple graph. For any vertex $v\in V$, the {\it
open neighborhood} of $v$ is the set $N_G(v)=\{u\in V | \{u, v\}\in E\}$
and the {\it closed neighborhood}
 is the set $N_G[v]=N_G(v)\cup \{v\}$.
 For a set $S\subseteq V$, the open neighborhood of $S$ is $N_G(S)=\bigcup_{v\in S} N_G(v)$ and the
  closed neighborhood of $S$ is $N_G[S]=N_G(S)\cup S$.
  For every vertex $v\in V (G)$, the {\it degree} of $v$ is the number of
edges incident with $v$ and is denoted by $d_G(v)$. We denote the maximum degree and the minimum degree of vertices of $G$ by $\Delta(G)$ and
$\delta(G)$, respectively.
A set $S\subseteq V$ is a {\it dominating set} if $N_G[S]=V$, or
equivalently,
 every vertex in $V\backslash S$ is adjacent to at least one vertex in $S$.
 The {\it domination number} $\gamma(G)$ is the minimum cardinality of a dominating set in $G$.
For a detailed treatment of these parameters, the reader is referred to~\cite{domination}.

 A dominating set with cardinality $\gamma(G)$ is called a {\it $\gamma$-set}.
 Let $D$ and $U$ be two subsets of $V$. The set $U$ is called a {\it monitor } set of $D$ if $D\subseteq N[U]$. The monitor number of $D$,
 denoted by $\iota(D)$, is the minimum cardinality of a monitor set of $D$. Set $\iota(G) = min \{\iota(D) :~D~ is~a~\gamma-set~ of ~G\}$.
 An {\it $i$-subset} of $V(G)$ is a subset of $V(G)$ with  cardinality $i$.
 Let ${\cal D}(G,i)$ be the family of
 dominating sets of $G$ which are $i$-subsets and let $d(G,i)=|{\cal D}(G,i)|$.
 The polynomial $D(G,x)=\sum_{i=0}^{|V(G)|} d(G,i) x^{i}$ is defined as {\it domination polynomial} of $G$ \cite{euro,saeid1}.
 Thus $D(G,x)$ is the generating polynomial for the number of dominating sets of $G$
of each cardinality.
Dominating sets are an important invariant of a graph and have many useful applications \cite{hook}.
   For more information and
motivation of domination polynomial   refer to
\cite{euro,saeid1}.

 Although the domination polynomial has been actively studied in recent
years, almost no attention has been given to the domination polynomials of some  graph
products. In fact the domination polynomials of binary graph operations, aside from union, join
and corona, have not been widely studied (see \cite{operation,peter}).

 Two graphs $G$ and $H$ are disjoint
if $V(G)\cap V(H)=\emptyset$. The union of two disjoint graphs $G$ and $H$ is the graph $G\cup H$
with the vertex set $V(G \cup H)= V(G) \cup V(H)$ and the edge set $E(G \cup H)=E(G) \cup E(H)$.
The notation $nG$ is the short notation for the union of $n$ copies of disjoint graphs isomorphic to $G$.

 The join of
two graphs $G_1$ and $G_2$, denoted by $G_1 \lor G_2$, is a graph with vertex set $V(G_1)\cup V(G_2)$ and edge set
$E(G_1)\cup  E(G_2)\cup \{uv |u\in V(G_1),  v \in V(G_2)\}$.

 For two graphs $G$ and $H$, let $G[H]$ be the graph with
vertex set $V(G)\times V(H)$,  such that the vertex $(a,x)$ is
adjacent to vertex $(b, y)$ if and only if $a$ is adjacent to $b$
(in $G$) or $a = b$ and $x$ is adjacent to $y$ (in $H$). The graph
$G[H]$ is the lexicographic product  of $G$ and
$H$.

 This product was introduced as the composition of graphs by Harary \cite{Har59,Har69}.
The lexicographic product is also known as graph
substitution, a name that bears witness to the fact that $G[H]$ can be obtained from $G$ by
substituting a copy $H_u$ of $H$ for every vertex $u$ of $G$ and then joining all vertices of $H_u$ with
all vertices of $H_v$ if $\{u,v\} \in E(G)$.

Calculating the domination number $\gamma (G)$ of the graph, and determining whether $\gamma (G) \leq k$ is known to be NP-complete \cite{garey}.
 Obviously, the problem is more complicated when we consider the lexicographic product of two graphs. The following result is about the domination number of lexicographic of two graphs.

\begin{figure}[!h]
\hspace{.05cm}
\includegraphics[width=13cm,height=4.1cm]{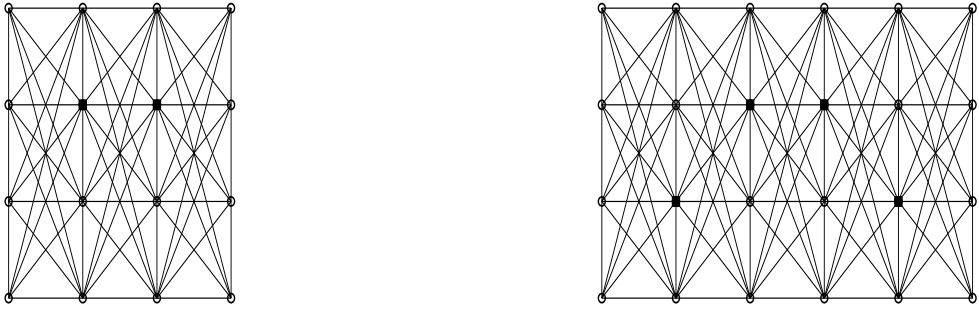}
\caption{ \label{fig}A dominating set of $P_4[P_4]$ and  $P_6[P_4]$, respectively. }
\end{figure}

\begin{theorem}\rm\cite{lex}\label{t}
\begin{enumerate}
\item[(i)]
Let $G$ and $H$ be two graphs with at least two vertices. If $\gamma(H)=1$, then $\gamma(G[H])=\gamma(G)$.
\item[(ii)]
Let $G$ be a graph with no isolated vertex, and let $H$ be a graph with $\gamma(H)\geq 2$. Then $\gamma(G)\leq \gamma(G[H]) \leq \gamma(G) +
\iota(G)$.
\end{enumerate}
\end{theorem}
 The lower bound and upper bound in Theorem \ref{t}(ii) are sharp.
Clearly, $\gamma(P_4)=2,~\gamma(P_6)=2,~\iota(P_4)=1,~\iota(P_6)=2.$ We have
$\gamma(P_4[P_4])=\gamma(P_4)=2,~\gamma(P_6[P_4])=\gamma(P_6)+\iota(P_6)=4.$ Thus, a domination number of $P_4[P_4]$ achieves the lower bound
and a domination number of $P_6[P_4]$ achieves the upper bound (see Figure \ref{fig}).

\medskip

 In this paper, we study  the domination polynomial of lexicographic product of specific graphs.

\section{Main results}

 In this section, first  we study some properties of lexicographic product of two graphs and then,
 we study the domination polynomial of lexicographic product of specific graphs.

 We mention  some properties of  lexicographical product of graphs.
 Note that the lexicographic product is not commutative, but it is associative and it is easily seen to have $K_1$ as both
a left and right unit. Also  the right-distributive rule,  holds for all
graphs $G, H$ and $K$:
\begin{eqnarray*}
(G_1[G_2])[G_3] &\cong &G_1([G_2[G_3]]),\\
K_1[G]& \cong &G,\\
G[K_1]&\cong &G,\\
(G \cup H)[K] &=& G[K] \cup H[K].
\end{eqnarray*}

 With respect to taking complements, we note that
$\overline{G[H]} = \bar{G}[\bar{H}]$ and since $\overline{\bar{G}} = G$, we have $G[H] = \overline{\bar{G}[\bar{H}]}$.
 We state and prove (with new approach) the following theorem which  has appeared  in \cite{Hamm} as Proposition 10.1.

\begin{theorem}
For a graph $G$ and natural number  $n\geq 2$,
\begin{enumerate}
\item[(ii)] $~G[nK_1]\cong (nK_1)[G]$ if and only if $G$ is totally disconnected,

\item[(i)] $ ~G[K_n]\cong K_n[G]$ if and only if $G$ is complete graph.
\end{enumerate}
\end{theorem}

\proof
\begin{enumerate}
\item[(i)] Observe that  $G[nK_1]\cong (nK_1)[G]$ is true, if $G$ has no edges.
Assume now that $G[nK_1]\cong (nK_1)[G]$. Evidently,
$$|E(G)| n^2 = |E(G[nK_1])| = |E(nK_1)[G]| = n |E(G)|.$$
So $|E(G)|(n^2 - n) = 0$. For $n\geq 2$, this is only possible if $|E(G)| = 0$.

\item[(ii)] We know that $\overline{G[H]} = \bar{G}[\bar{H}].$ Therefore $G$ and $K_n$ commute if and only if
$\bar{G}$ and $\bar{K_n}$ commute. For $n\geq 2$ this is true  if and only if $\overline{G}$ has no edges, or
equivalently if and only if $G$ is complete graph.\qed
\end{enumerate}

\begin{lemma}\label{lem1}
 For every graph $G$ and  natural numbers $m$ and $n$, we have:
\begin{enumerate}
\item[$(i)$] $(nK_1)[G] \cong nG$ and $K_m[K_n] \cong K_{mn}$.
\item[$(ii)$]  $K_m[G] \cong \lor _{i=1}^m G$. In particular, if $m=2$ and $G=nK_1$, then $K_2[nK_1] \cong K_{n,n}$.
\end{enumerate}
\end{lemma}

\proof
\begin{enumerate}
\item[(i)] It follows from  the definition of lexicographic product of two graphs.

\item[(ii)] By the definition of lexicographic product, the graph $K_m[G]$  is the graph obtained from $mG$ such that  all vertices of copies of $G$ are
    adjacent. 
Since  $K_{n,n} \cong (nK_1)\lor(nK_1)$, we have $K_2[nK_1] \cong K_{n,n}$.\qed
 \end{enumerate}

\begin{figure}[h]
\hspace{3.cm}
\includegraphics[width=8.5cm,height=2.cm]{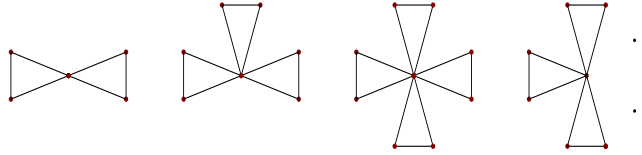}
\caption{\label{figure1} Friendship graphs $F_2, F_3, F_4$ and $F_n$, respectively.}
\end{figure}

Here, we recall the definition of friendship graphs.  The friendship (or Dutch-Windmill) graph $F_n$ is a graph that can be constructed by coalescence $n$
copies of the cycle graph $C_3$ of length $3$ with a common vertex \cite{site}. See Figure \ref{figure1}.
Domination polynomials, exploring the nature and location of roots of domination polynomials of friendship graphs has studied in \cite{jason}.

 For any two graphs $G$ and $H$, if $G$ is an isolated vertex, then $G[H]\cong H$, and if $H$ is an isolated vertex, then $G[H]\cong  G$.
  Note that the right-distributive rule $(K \lor H)[G] = K[G] \lor H[G]$ holds for all
graphs $G, H$ and $K$. In particular, $(K_1 \lor H)[G] = G \lor H[G]$. Now, we state and prove the following theorem which is about the lexicographic product of star and friendship graph with an arbitrary graph.
 \begin{lemma}\label{lem3}
For every graph $G$ and  natural number $n$,
\begin{enumerate}
\item[(i)] $K_{1,n}[G] \cong G \lor nG$.

\item[(ii)]  $F_{n}[G] \cong G \lor n(G\lor G)$.
\end{enumerate}
\end{lemma}

\proof
\begin{enumerate}
\item[(i)] Let $u$ be a center of $K_{1,n}$. By the definition of lexicographic product, the graph   $K_{1,n}[G]$  is the graph  $(n+1)G$ together with all edges  that joining all vertices of $G_u$ with all vertices of each other copy of $G$.
In other words,  $K_{1,n}[G] \cong G \lor nG$.
\item[(ii)] Since the friendship graph  $F_n$ is the join of $K_1$ and $nK_2$,  we have the result.\qed
 \end{enumerate}


 To obtain the domination polynomials of lexicographic product of some graphs,  we need  some
preliminary properties of the  domination polynomial.

\begin{theorem}\label{theorem4}{\rm \cite{ saeid1}}
If a graph $G$ consists of $k$ components $G_1, \dots ,G_k$,  then
$D(G, x) =\prod_{i=1}^k D(G_i, x).$
\end{theorem}

The following theorem gives the domination polynomial of join of some graphs.

\begin{theorem}\label{join}{\rm\cite{euro}}
Let $G_1, \ldots, G_k$ be  graphs of orders $n_1, \ldots, n_k$,
respectively. Then
\[
D(\lor_{i=1}^{k} G_i,x)=\sum_{j=1}^{k-1}\Big((1+x)^{n_j}-1\Big) \Big((1+x)^{\sum_{i=j+1}^{k}n_i}-1\Big)+ \sum_{i=1}^kD(G_i,x).
\]
\end{theorem}

 The following results gives the  domination polynomials of lexicographic product with complete graphs.

\begin{theorem}\label{lem2}
For every graph $G$ and every natural number $m$ and $n$,
\begin{enumerate}
\item[$(i)$] $D((nK_1)[G],x) = (D(G,x))^n$.
\item[$(ii)$] $D(K_m[K_n],x) = (1+x)^{mn}-1$.
\item[$(iii)$] $D(K_m[G],x) = \Big((1+x)^{|V(G)|}-1\Big)\sum_{j=1}^{m-1} \Big((1+x)^{(m-j)|V(G)|}-1\Big)+ mD(G,x)$.
\end{enumerate}
\end{theorem}

\proof
\begin{enumerate}
\item[(i)] By Lemma \ref{lem1}(i), $D((nK1)[G], x) = D((nK1),D(G,x))$. Since $D(nK_1,x)=x^n$, we have the result.

\item[(ii)] It is easy to see that by Lemma \ref{lem1}(i), $D(K_m[K_n],x) = D(K_m,D(K_n,x))$. Since $D(K_t,x)=(1+x)^t-1$, we have the result.

\item[(iii)] Using   Lemma \ref{lem1}(iii) and Theorem \ref{join}, we have
\begin{eqnarray*}
D(\lor_{i=1}^{m} G,x)&=&\sum_{j=1}^{m-1}\Big((1+x)^{|V(G)|}-1\Big) \Big((1+x)^{\sum_{i=j+1}^{m}|V(G)|}-1\Big)+ \sum_{i=1}^mD(G,x)\\
&=&\Big((1+x)^{|V(G)|}-1\Big)\sum_{j=1}^{m-1} \Big((1+x)^{(m-j)|V(G)|}-1\Big)+ mD(G,x).
\end{eqnarray*}
\end{enumerate}

Using  Lemma \ref{lem3},  Theorems \ref{theorem4} and  \ref{join}, we have the following theorem which  is
 about the domination polynomial of lexicographic product with star and friendship  graphs.

\begin{theorem}
\begin{enumerate}

\item[(i)]

$D(K_{1,n}[G],x) = \Big((1+x)^{|V(G)|}-1\Big)\Big((1+x)^{n|V(G)|}-1\Big)+D(G,x)+D(G,x)^{n}.$
\item[(ii)]
\begin{eqnarray*}
D(F_{n}[G],x)&=&\Big((1+x)^{|V(G)|}-1\Big)\Big((1+x)^{2n|V(G)|}-1\Big)+D(G,x)\\
&&+\Big(\big((1+x)^{|V(G)|}-1\big)^2+2D(G,x)\Big)^{n}.
\end{eqnarray*}
\end{enumerate}

\end{theorem}

  For some times, we thought that $D(G[H],x)=D(G,D(H,x)-1)$ is true \cite{ENDM} (similar to independence domination polynomial of a graph
\cite{Hickman}), but Theorem \ref{lem2} show that this formula is not true.  Note that in \cite{operation,brown} it has proved that $D(G[K_n], x)=D(G, D(K_n,x))$, which is an extension of Theorem \ref{lem2}(ii).

\medskip

\section{Conclusion}

In this paper, we studied the domination polynomial of lexicographic product of two specific graphs, but for two arbitrary graphs the problem
is open yet. Until now all attempts to find a  formula for $D(G[H],x)$  failed.

\end{document}